\documentclass{article}%
\usepackage{amssymb}
\usepackage{amsmath}
\usepackage{amsfonts}
\usepackage{amsthm}
\usepackage{graphicx}%
\providecommand{\U}[1]{\protect\rule{.1in}{.1in}}
\newtheorem{theorem}{Theorem}[section]

\newtheorem{example}[theorem]{Example}
\begin{document}

\title{On the primitive element of finite $\Bbbk-$ algebras and applications to
commuting matrices}
\author{Aristides I. Kechriniotis \thanks{Department of Physics, University of
Thessaly, GR-35100 Lamia, Greece. Corresponding author. E-mail:
arisk7@gmail.com }}
\date{March 18, 2024}
\maketitle

\begin{abstract}
Let $\Bbbk$ be a field of characteristic zero. Using the properties of the
ideal of the coordinate Hermite interpolation on an $n$-dimensional grid [4],
we prove that the extension $\Bbbk\subset\Bbbk\left[  x_{1},x_{2}%
,...,x_{n}\right]  /\left(  f_{1}\left(  x_{1}\right)  ,...,f_{n}\left(
x_{n}\right)  \right)  $ has a primitive element if and only if at most one of
the univariate polynomials $f_{1},...,f_{n}$ is inseparable. This result
\ leads to some interesting conclusions regarding the primitive elements of
finite $\Bbbk-$ algebras. Finally, these results are further used to
investigate the Frobenius' question, of whether two commuting matrices $A$ and
$B~$can be expressed as polynomials in some matrix $C.$

\end{abstract}


\section{Introduction}%

\numberwithin{equation}{section}%
The main result of this note is contained within the following Theorem:

\begin{description}
\item[Theorem 1.1] \textit{Let }$\Bbbk$\textit{\ be a field of characteristic
zero, and }$I$\textit{\ an ideal of }$\Bbbk\left[  x_{1},...,x_{n}\right]  $,
$n\geq2$\textit{\ generated by the univariate polynomials }$f_{1}\in
\Bbbk\left[  x_{1}\right]  ,...,f_{n}\in\Bbbk\left[  x_{n}\right]  $\textit{.
Then the finite }$\Bbbk-$\textit{\ algebra }$\Bbbk\left[  x_{1},...,x_{n}%
\right]  /I$\textit{\ has a primitive element over}$\mathit{\ }\Bbbk$,
\textit{if and only if at most one of the polynomials }$f_{1},...,f_{n}%
$\textit{\ is inseparable.}
\end{description}

\bigskip By Theorem 1.1 we easily get the following, related to the existence
of primitive elements of finite $\Bbbk-$algebras, Corollaries :

\begin{description}
\item[Corollary 1.2] \textit{Let }$\Delta$\textit{\ be a finite }$\Bbbk
-$\textit{algebra. Given }$\delta_{1},...,\delta_{n}\in\Delta$\textit{\ such
that at most one of their minimal polynomials }$\mu_{\delta_{1}}%
,...,\mu_{\delta_{n}}$\textit{\ is inseparable, then the extension }%
$\Bbbk\subset\Bbbk\left[  \delta_{1},...,\delta_{n}\right]  $\textit{\ has a
primitive element.}

\item[Corollary 1.3] \textit{Let }$\Delta$\textit{\ be as above. Given
}$\delta_{1},\delta_{2}\in\Delta$\textit{\ such that }$\mu_{\delta_{1}}%
,\mu_{\delta_{2}}$\textit{\ are inseparable, and }%
\begin{equation}
\deg\mu_{\delta_{1}}\deg\mu_{\delta_{2}}=\dim_{\Bbbk}\Bbbk\left[  \delta
_{1},\delta_{2}\right]  , \tag{1.1}%
\end{equation}
\textit{then the extension }$\Bbbk\subset\Bbbk\left[  \delta_{1},\delta
_{2}\right]  $\textit{\ has no any primitive element.}
\end{description}

Frobenius [1] posed the following question: Can any two commuting matrices $A$
and $B$ over a field $\Bbbk$ be expressed as polynomials in some matrix $C$?

By considering the following two commuting matrices
\begin{equation}
A=\left[
\begin{array}
[c]{ccc}%
0 & 1 & 0\\
0 & 0 & 0\\
0 & 0 & 0
\end{array}
\right]  \text{, }B=\left[
\begin{array}
[c]{ccc}%
0 & 0 & 1\\
0 & 0 & 0\\
0 & 0 & 0
\end{array}
\right]  , \tag{1.2}%
\end{equation}
\ he was able to prove that there is no matrix $C$ $\in\Bbbk\left[
A,~B\right]  $, such that $A$ and $B$ can be expressed as polynomials in $C$,
and posed the following question: Is there some matrix $C$, such that any pair
of commuting matrices $A$ and $B$ over a field $\Bbbk$ can be expressed as
polynomials in $C$?

Shoda [5] has shown that the matrices $A$ and $B$ as given in $\left(
1.2\right)  $ can not be expressed as polynomials in any matrix $C$. In the
example 3.1 in section 3, we also provide a very simple proof that the above
$3\times3$ matrices $A$ and $B$ cannot be written as polynomials in any matrix
$C.$

Two main results related to the Frobenius question can be found in the literature:

1. If $A$ commutes with $B$, and the minimal polynomial of $A$ is equal to its
characteristic polynomial, then $B$ can be expressed as polynomial in $A$,

[2], [3; pp.135-137 ].

2. If two commuting matrices $A$ and $B$ are both diagonalizable, then $A$ and
$B$ can be expressed as polynomials in some matrix $C$.

In the present paper, using Corollaries 1.2 and 1.3 the following results are obtained:

Let $\Bbbk$ be a field of characteristic zero, and two $n\times n$ matrices
$A$ and $B$ over the field $\Bbbk.$

\begin{description}
\item[Theorem 1.4] \textit{Let }$A_{1},...,A_{m},m\geq2$\textit{\ be commuting
}$n\times n$\textit{\ matrices over }$\Bbbk$\textit{. If at most one of the
matrices }$A_{1},...,A_{m}$\textit{\ is not diagonalizable, then all matrices
can be expressed as polynomials in some matrix }$C.$
\end{description}

Theorem 1.4 for $m=2$ improves the above result 2., and may be formulated as
follows: If at least one of the commuting matrices $A\ $\ and $B$ is
diagonalizable, then $A$ and $B$ can be expressed as polynomials in some
matrix $C$.

\begin{description}
\item[Theorem 1.5] Let $n$ be any composite number. \textit{If the minimal
polynomials }$\mu_{A}$\textit{\ and }$~\mu_{B}$\textit{\ of \ the commuting
}$n\times n$\textit{ matrices }$A$\textit{\ and }$B$\textit{\ are both
inseparable, and }$\deg\mu_{A}\deg\mu_{B}=\dim_{\Bbbk}\Bbbk\left[  A,B\right]
$\textit{, then neither }$A$\textit{\ nor }$B$\textit{\ can be expressed as
polynomials in some matrix }$C\in\Bbbk\left[  A,B\right]  .$ In addition, if
\begin{equation}
\deg\mu_{A}\deg\mu_{B}=\dim_{\Bbbk}\Bbbk\left[  A,B\right]  =n, \tag{1.3}%
\end{equation}
\textit{\ then neither }$A$\textit{\ nor }$B$\textit{\ can be expressed as
polynomials in some matrix }$C.$
\end{description}

According to Corollary 3.2, we can construct the set of all pairs of commuting
and inseparable matrices satisfying $(1.3).$

\section{The primitive element Theorem}

\bigskip Let $\Bbbk$ be $~$a field of characterstic zero and $\overline{\Bbbk
}$ the algebraic closure of $\Bbbk.$For convenience we define the following
symbols and notations:

\begin{itemize}
\item $%
\mathbb{N}
_{0}:=\left\{  0,1,2,...\right\}  .$

\item $\left\vert A\right\vert $ is the cardinality of the set $A$;

\item $A^{n}:=\underset{n-times}{A\times...\times A}$,

\item Given the sets $A_{1},...,A_{n}$, then $\mathbf{A:=}A_{1}\times...\times
A_{n}$. Further any element $\left(  c_{1},...c_{n}\right)  \in\mathbf{A}$
will be denoted by $\mathbf{c.}$

\item $\mathbf{0:=}\left(  0,0,...,0\right)  $, \ $\mathbf{1:=}\left(
1,1,...,1\right)  .$

\item For $\mathbf{a\in}\Bbbk^{n}$, and $\mathbf{m\in}%
\mathbb{N}
_{0}^{n}$ we denote $\mathbf{a}^{\mathbf{m}}:=\prod_{i=1}^{n}a_{i}^{m_{i}}.$

\item In $%
\mathbb{N}
_{0}^{n}$ we define the relation "$\leq"$ as follows: $\mathbf{k}%
\leq\mathbf{m}$ if and only if $k_{i}\leq m_{i}$, for every $i=1,...,n.$
Clearly $\left(
\mathbb{N}
_{0}^{n},\leq\right)  $ is a poset ($%
\mathbb{N}
_{0}^{n}$ is partially ordered).

\item If $\mathbf{k\leq m}$ and $\mathbf{k\neq m}$, then $\left[
\mathbf{k,m}\right]  :=\left\{  \mathbf{l\in}%
\mathbb{N}
_{0}^{n}:\mathbf{k\leq l\leq m}\right\}  =\left[  k_{1},m_{1}\right]
\times...\times\left[  k_{n},m_{n}\right]  $, and is valid $\left\vert \left[
\mathbf{k,m}\right]  \right\vert =\prod_{i=1}^{n}\left(  m_{i}-k_{i}+1\right)
$ .

\item $\partial_{i}^{k}:=\frac{\partial^{k}}{\partial x_{i}^{k}}$,
$\partial^{\mathbf{k}}:=\prod_{i=1}^{n}\partial_{i}^{k_{i}}$, $\partial
_{\mathbf{a}}^{\mathbf{k}}f\left(  \mathbf{x}\right)  :=\partial^{\mathbf{k}%
}f\left(  \mathbf{x}\right)  \left\vert _{\mathbf{x=a}}\right.  .$

\item Let $\Delta$ be a finite $\Bbbk-$algebra. The minimal polynomial of
$\delta\in\Delta$ over $\Bbbk$ will be denoted by $\mu_{\delta}$, and the
number $\max\left\{  \deg\mu_{\delta}:\delta\in\Delta\right\}  $ will be
denoted by $\operatorname{codim}_{\Bbbk}\Delta.$
\end{itemize}

\begin{description}
\item[Remark 2.1] \textit{Let }$\Delta$\textit{ be as above. It is easy to
verify, that} \textit{the extension }$\Bbbk\subset\Delta$\textit{\ has a
primitive element if and only if }$\operatorname{codim}_{\Bbbk}\Delta
=\dim_{\Bbbk}\Delta$\textit{. This will be used in the proof of Theorem 1.1}.
\end{description}

To prove the main Theorem 1.1 some Lemmas are also required.

\begin{description}
\item[Lemma 2.2] \textit{Let }$A_{1},...,A_{n}$\textit{\ be finite subsets of
}$\overline{\Bbbk}$\textit{\ such that }$\left\vert A_{i}\right\vert
>1,i=1,...,n$\textit{. Then there exist non zero }$c_{2},...,c_{n}\in\Bbbk
$\textit{\ , such that the restriction of }$g:=x_{1}+c_{2}x_{2}+...+c_{n}%
x_{n}\in\Bbbk\left[  x_{1},...,x_{n}\right]  $\textit{\ on }$A$\textit{\ is an
injective mapping.}
\end{description}

\begin{proof}
We define in $\overline{\Bbbk}\left[  x_{1},...,x_{n}\right]  $ the polynomial
\[
f\left(  \mathbf{x}\right)  :=\prod_{_{\substack{\mathbf{a},\mathbf{b}%
\in\mathbf{A}\\\mathbf{a}\neq\mathbf{b}}}}\sum_{i=1}^{n}\left(  a_{i}%
-b_{i}\right)  x_{i}.
\]
Clearly $f$ is not identically zero. Suppose $\Bbbk^{n}\cap
\operatorname*{supp}\left(  f\right)  =\varnothing.$ Note that all
coefficients of the polynomial $f$ are algebraic over $\Bbbk$. Let $M$ be an
extension field of $\Bbbk$, generated from all coefficients of the polynomial
$f.$ Then $\Bbbk\subset M$ is a finite extension. Let $\dim_{\Bbbk}M=s$. Then
the extension $\Bbbk\subset M$ has a primitive element $\rho\in M$ . Also
$1,\rho,...,\rho^{s-1}$are linearly independent over $\Bbbk$. Therefore there
are unique polynomials $p_{i}\in\Bbbk\left[  x_{1},...,x_{n}\right]
,~i=0,1,...,s-1$ such that $f=\sum_{i=0}^{s-1}p_{i}\rho^{i}$. From $\Bbbk
^{n}\cap\operatorname*{supp}\left(  f\right)  =\varnothing$ it follows that
$\sum_{i=0}^{s-1}p_{i}\rho^{i}=0$ on $\Bbbk^{n}$. Therefore all polynomials
$p_{i},~i=0,1,...,s-1$ are identically zero. Thus we have that $f=0.$ This
contradics to $f$ is not identically zero. Consequently there is at least
one $\mathbf{d}=\left(  d_{1},...,d_{n}\right)  \in\Bbbk^{n}$ , such that
$f\left(  \mathbf{d}\right)  \neq0.$ Consider the linear mapping $h$ defined
by $\ h\left(  \mathbf{x}\right)  :=\sum_{i=1}^{n}d_{i}x_{i}$. Then from
$f\left(  \mathbf{d}\right)  \neq0$ it follows that%
\[
\prod_{\substack{\mathbf{a},\mathbf{b}\in\mathbf{A}\\\mathbf{a}\neq\mathbf{b}%
}}\left(  h\left(  \mathbf{a}\right)  -h\left(  \mathbf{b}\right)  \right)
\neq0.
\]
That means that the restriction of $\ h$ on $\mathbf{A}$ is an injective
mapping. Further, if we assume that $d_{i}=0$ for some $i\in\left\{
1,...,n\right\}  $, then $h$ is not an injective mapping on $\mathbf{A}.$
Therefore the following holds: $\prod_{i=1}^{n}d_{i}\neq0$ . Thus we have
\[
\operatorname*{card}\left(  d_{1}^{-1}h\left(  \mathbf{A}\right)  \right)
=\operatorname*{card}\left(  h\left(  \mathbf{A}\right)  \right)
\]%
\[
=\operatorname*{card}\left(  \mathbf{A}\right)  .
\]
Therefore $g:$ $\mathbf{A}\rightarrow\overline{\Bbbk},$ $\ g\left(
\mathbf{x}\right)  :=x_{1}+\sum_{i=1}^{n}c_{i}x_{i}$, where $c_{i}:=d_{i}%
d_{1}^{-1},$ is also an injective mapping on $\mathbf{A}.$
\end{proof}

\begin{description}
\item[Lemma 2.3] \textit{Given the real numbers }$k_{1},...,k_{n}%
,$\textit{\ such that }$k_{1}\geq1,...,k_{n}\geq1,$ \textit{we have that
}$\sum_{i=1}^{n}k_{i}-\left(  n-1\right)  =\prod_{i=1}^{n}k_{i}$\textit{\ if
and only if at most one number of }$k_{1},...,k_{n}$\textit{\ is greater than
}$1$\textit{.}
\end{description}

\begin{proof}
If $n-1$ numbers from $k_{1},...,k_{n}$ are equal to $1$, without loss of
generallity we can assume that $k_{1}=k\geq1,$ and $k_{2}=k_{3}=...=k_{n}=1$.
In this case, we have that $\sum_{i=1}^{n}k_{i}-\left(  n-1\right)  =k=\prod_{i=1}^{n}%
k_{i}.$ Therefore, if at most one number from $k_{1},...,k_{n}$\ is greater than
$1,$  we have that $\sum_{i=1}^{n}k_{i}-\left(  n-1\right)  =\prod_{i=1}%
^{n}k_{i}.$ Furthermore, if we suppose that $\ r>1$ numbers from $k_{1},...,k_{n}$ are
greater than $1,$ without loss of generality we can assume that $k_{1}%
>1,...,k_{r}>1$ and $k_{r+1}=k_{r+2}=...=k_{n}=1.$ In this case, we obtain that
\begin{align*}
\sum_{i=1}^{n}k_{i}-\left(  n-1\right)   &  =\sum_{i=1}^{r}k_{i}-r+1\\
&  =\sum_{i=1}^{r}\left(  k_{i}-1\right)  +1\\
&  <\prod_{i=1}^{r}\left(  1+\left(  k_{i}-1\right)  \right)  \\
&  =\prod_{i=1}^{r}k_{i}=\prod_{i=1}^{n}k_{i},
\end{align*}
and consequently $\sum_{i=1}^{n}k_{i}-\left(  n-1\right)  \neq\prod_{i=1}^{n}k_{i}.$
Therefore, at most one number from $k_{1},...,k_{n}$ satisfying the condition $\sum
_{i=1}^{n}k_{i}-\left(  n-1\right)  =\prod_{i=1}^{n}k_{i}$\ could be greater
than $1.$ Summarizing the above we immediately get the claim of Lemma 2.3.
\end{proof}

Given the univariate polynomials $f_{i}\in\Bbbk\left[  x_{i}\right]
,~i=1,...,n\,,$ we denote by $A_{i}$ the set of the roots of $f_{i}$ in the
algebraic closure $\overline{\Bbbk}$ of $\Bbbk$ and by $\nu_{i}\left(
a\right)  $ the multiplicity of $a\in A_{i}.$That is $f_{i}\left(
x_{i}\right)  =\prod_{a\in A_{i}}\left(  x_{i}-a\right)  ^{\nu_{i}\left(
a\right)  },~i=1,...,n;.$

We denote by $I$ the ideal of $\Bbbk\left[  x_{1},...,x_{n}\right]  $
generated by $f_{1},...,f_{n}$, and by $\overline{I}$ the ideal of
$\overline{\Bbbk}\left[  x_{1},...,x_{n}\right]  $ generated also by
$f_{1},...,f_{n}.$

\bigskip For any $g\in\Bbbk\left[  x_{1},...,x_{n}\right]  $ let us denote
by\ $\widetilde{g}$ the residue class of $g$ in $\Bbbk\left[  x_{1}%
,...,x_{n}\right]  /I$ , and by $\overline{g}\ $the residue class of $g$\ in
$\overline{\Bbbk}\left[  x_{1},...,x_{n}\right]  /\overline{I}$ .

From [4] we need some results included in the following two Lemmas.

\begin{description}
\item[Lemma 2.4] \textit{(see [4], Lemma 4.2) Let }$I=\left(  f_{1}%
,...,f_{n}\right)  $\textit{\ be\ an ideal of }$\Bbbk\left[  x_{1}%
,...,x_{n}\right]  $\textit{\ generated by the polynomials }$f_{1}\in
\Bbbk\left[  x_{1}\right]  ,...,f_{n}\in\Bbbk\left[  x_{n}\right]
.$\textit{\ Then}
\begin{align*}
\dim_{\Bbbk}\Bbbk\left[  x_{1},...,x_{n}\right]  /I  &  =\prod_{i=1}^{n}\deg
f_{i}=\prod_{i=1}^{n}\sum_{a_{i}\in A_{i}}\nu_{i}\left(  a\right) \\
&  =\dim_{\overline{\Bbbk}}\overline{\Bbbk}\left[  x_{1},...,x_{n}\right]
/I\overline{\Bbbk}\left[  x_{1},...,x_{n}\right]  .
\end{align*}

\item[Lemma 2.5] \textit{(see [4], Corollary 4.4) }$f\in I$\textit{\ if and
only if }$\ \partial_{\mathbf{a}}^{\mathbf{m}}f=0$\textit{\ for every
}$\mathbf{a}\in\mathbf{A}$\textit{\ and}$~\mathbf{m}\in\left[  \mathbf{0,}%
\nu\left(  \mathbf{a}\right)  -\mathbf{1}\right]  $, where $\nu\left(
\mathbf{a}\right)  :=\left(  \nu_{1}\left(  a_{1}\right)  ,...,\nu_{n}\left(
a_{n}\right)  \right)  .$

\item[Lemma 2.6] \textit{Given }$g\in\Bbbk\left[  x_{1},...,x_{n}\right]
$\textit{, the following expressions are equivalent:}

\item[1] \textit{For any }$f\in\Bbbk\left[  x_{1},...,x_{n}\right]
$\textit{\ there is one }$p\in\Bbbk\left[  t\right]  $\textit{\ such that
}$f=p\left(  g\right)  +I$\textit{.}

\item[2] \textit{For any }$h\in\overline{\Bbbk}\left[  x_{1},...,x_{n}\right]
$\textit{\ there is one }$q\in\overline{\Bbbk}\left[  t\right]  $%
\textit{\ such that }$h=q\left(  g\right)  +\overline{I}.$
\end{description}

\begin{proof}
Let the first expression be true . Let $L$ be the extension field of $\Bbbk$
generated from all coefficients of $h\,.$ Clearly the algebraic extenssion
$\Bbbk\subset L~$ is finite. Let $\dim_{\Bbbk}L=r.$ Then the field extenssion
$\Bbbk\subset L$ has a primitive element $\rho\in L$. Thus each coefficient
$c$ of $h~$can be written as $c=\sum_{i=0}^{r-1}a_{i}\rho^{i},$ where
$\alpha_{1},...,a_{n-1}\in\Bbbk$. Therefore $h$ takes the following form:$~$%
\begin{equation}
h=\sum_{i=0}^{r-1}f_{i}\rho^{i},\tag{2.1}%
\end{equation}
where $f_{i}\in\Bbbk\left[  x_{i}\right]  \subset\Bbbk\left[  x_{1}%
,...,x_{n}\right]  .$ From expression 1, there are $p_{i}\in\Bbbk\left[
t\right]  ,$ such that \ for every $~i=1,...,n$
\begin{equation}
f_{i}=p_{i}\left(  g\right)  +I.\tag{2.2}%
\end{equation}
Substituting $\left(  2.2\right)  $ into $\left(  2.1\right)  $ we get the second expression.
Let the second expression be true. Then there are $q\in\overline{\Bbbk}\left[
t\right]  ,$ and $~h_{1},...,h_{n}\in\overline{\Bbbk}\left[  x_{1}%
,...,x_{n}\right]  $ such that
\begin{equation}
f=q\left(  g\right)  +\sum_{i=1}^{n}h_{i}f_{i}.\tag{2.3}%
\end{equation}
Let $M$ be an extenssion field of $\Bbbk$ generated from all coefficients of
$q,h_{1},...,h_{n}$. Then $\Bbbk\subset M$ is finite. If $\dim_{\Bbbk}M=s$,
then $\Bbbk\subset M$ has a primitive element $\sigma\in M$. Working in a
similar manner as in the first part, $\left(  2.3\right)  $ can be rewritten
as:
\begin{align}
f &  =q_{0}\left(  g\right)  +\sum_{i=1}^{n}h_{i0}f_{i}+\left(  q_{1}\left(
g\right)  +\sum_{i=1}^{n}h_{i1}f_{i}\right)  \sigma\tag{2.4}\\
&  +...+\left(  q_{s-1}\left(  g\right)  +\sum_{i=1}^{n}h_{is-1}f_{i}\right)
\sigma^{s-1},\nonumber
\end{align}
where $q_{i}\in\Bbbk\left[  t\right]  $ for all $i=1,...,n$,~and $h_{ij}%
\in\Bbbk\left[  x_{1},...,x_{n}\right]  $ for all $i=0,1,...,n$,
$~j=1,...,s-1.$ Since $q_{0}\left(  g\right)  +\sum_{i=1}^{n}h_{i0}f_{i}%
-f,~$\ and $q_{j}\left(  g\right)  +\sum_{i=1}^{n}q_{ij}f_{i}~\in\Bbbk\left[
x_{1},...,x_{n}\right]  $, $~j=1,...,s-1$, and $\left\{  1,\sigma
,...,\sigma^{s-1}\right\}  $ is linearly indepedent over $\Bbbk$, from
$\left(  2.4\right)  ~$we get that $f=q_{0}\left(  g\right)  +\sum_{i=1}%
^{n}q_{i0}f_{i}.$ That is equivalent to $f=q_{0}\left(  g\right)  +I.$
\end{proof}

\begin{description}
\item[Lemma 2.7] \textit{For any }$f\in\overline{\Bbbk}\left[  x_{1}%
,...,x_{n}\right]  ,~$\textit{we have that }%
\[
F:=\prod_{\mathbf{a\in A}}\left(  f\left(  \mathbf{x}\right)  \mathbf{-f}%
\left(  \mathbf{a}\right)  \right)  ^{\sum_{i=1}^{n}\nu_{i}\left(
a_{i}\right)  -n+1}\in\overline{I}.
\]

\end{description}

\begin{proof}
Fix $\mathbf{a\in A.}$ Taking into account that there are $g_{i}\in
\overline{\Bbbk}\left[  x_{1},...,x_{n}\right]  $ such that $f\left(
\mathbf{x}\right)  \mathbf{-f}\left(  \mathbf{a}\right)  =\sum_{i=1}%
^{n}\left(  x_{i}-a_{i}\right)  g_{i}\left(  \mathbf{x}\right)  $, $F$ can be
written as%
\begin{equation}
F=g\left(  \mathbf{x}\right)  \left(  \sum_{i=1}^{n}\left(  x_{i}%
-a_{i}\right)  g_{i}\left(  \mathbf{x}\right)  \right)  ^{\sum_{i=1}^{n}%
\nu_{i}\left(  a_{i}\right)  -n+1},\tag{2.5}%
\end{equation}
for some $g\in\overline{\Bbbk}\left[  x_{1},...,x_{n}\right]  .$ Using the
multinomial formula on the right-hand side of $\left(  2.5\right)  $ we have$,$
\[
F=\sum_{\left\vert \mathbf{i}\right\vert =\sum_{i=1}^{n}\nu_{i}\left(
a_{i}\right)  -n+1}\left(  \mathbf{x}-\mathbf{a}\right)  ^{\mathbf{i}%
}h_{\mathbf{i}}\left(  \mathbf{x}\right)  ,
\]
for some $h_{\mathbf{i}}\in\overline{\Bbbk}\left[  x_{1},...,x_{n}\right]  .$
Note that for any $\mathbf{m\in}\left[  \mathbf{0,}\nu\left(  \mathbf{a}\right)
-\mathbf{1}\right]  $ the following is valid: $\left\vert \mathbf{m}%
\right\vert \leq\sum_{i=1}^{n}\nu_{i}\left(  a_{i}\right)  -n<\sum_{i=1}%
^{n}\nu_{i}\left(  a_{i}\right)  -n+1=\left\vert \mathbf{i}\right\vert .$
Therefore, for every $\mathbf{m\in}\left[  \mathbf{0,}\nu\left(  \mathbf{a}%
\right)  -\mathbf{1}\right]  $ and every $\mathbf{i}$ satisfying $\left\vert
\mathbf{i}\right\vert =\sum_{i=1}^{n}\nu_{i}\left(  a_{i}\right)  -n+1$ there
is at least one $j\in\left\{  1,...,n\right\}  $ such that $i_{j}>m_{j}.$ Now
applying the Leibniz derivative rule for multivariable functions we get that
$\partial_{\mathbf{a}}^{\mathbf{m}}F=0,~~\mathbf{m\in}\left[  \mathbf{0,}%
\nu\left(  \mathbf{a}\right)  -\mathbf{1}\right]  $, and by Lemma 2.5, it
follows that $F\in\overline{I}.$
\end{proof}

The main result of this paper can now be proven.

\begin{proof}
[Proof of Theorem 1.1]Let $A_{1},...,A_{n}$ be the sets of the zeros of
$f_{1},...,f_{n}$ respectively, and let $\nu_{i}\left(  a\right)  $ be the
multiplicity of $a\in A_{i}$. First it will be proved that
\begin{equation}
\operatorname{codim}_{\overline{\Bbbk}}\left(  \overline{\Bbbk}\left[
x_{1},...,x_{n}\right]  /\overline{I}\right)  =\sum_{\mathbf{a\in A}}%
\sum_{i=1}^{n}\left(  \nu_{i}\left(  a_{i}\right)  -n+1\right)  :\tag{2.6}%
\end{equation}
Let $f$ be any polynomial in $\overline{\Bbbk}\left[  x_{1},...,x_{n}\right]
$. Then by Lemma 2.7  follows that $\prod_{\mathbf{a\in A}}\left(  f\left(
\mathbf{x}\right)  \mathbf{-f}\left(  \mathbf{a}\right)  \right)  ^{\sum
_{i=1}^{n}\nu_{i}\left(  a_{i}\right)  -n+1}\in\overline{I}.$ Therefore ~
\[
\prod_{\mathbf{a\in A}}\left(  \overline{f\left(  \mathbf{x}\right)
}\mathbf{-f}\left(  \mathbf{a}\right)  \right)  ^{\sum_{i=1}^{n}\nu_{i}\left(
a_{i}\right)  -n+1}=\overline{0},
\]
where $\overline{f\left(  \mathbf{x}\right)  }$ \ $\in\overline{\Bbbk}\left[
x_{1},...,x_{n}\right]  /\overline{I}$ is the equivalence class of $f\left(
\mathbf{x}\right)  $. Therefore, there is a$~H\in\overline{\Bbbk}\left[
t\right]  $ such that%
\begin{equation}
\prod_{\mathbf{a\in A}}\left(  t\mathbf{-f}\left(  \mathbf{a}\right)  \right)
^{\sum_{i=1}^{n}\nu_{i}\left(  a_{i}\right)  -n+1}=\mu_{\overline{f\left(
\mathbf{x}\right)  }}\left(  t\right)  H\left(  t\right)  .\tag{2.7}%
\end{equation}
From $\left(  2.7\right)  ~$we conclude that $\deg\mu_{\overline{f\left(
\mathbf{x}\right)  }}\leq\sum_{\mathbf{a\in A}}\left(  \sum_{i=1}^{n}\nu
_{i}\left(  a_{i}\right)  -n+1\right)  $, for any $f\in$ $\overline{\Bbbk
}\left[  x_{1},...,x_{n}\right]  $. Therefore we have that
\begin{equation}
\operatorname{codim}_{\overline{\Bbbk}}\left(  \overline{\Bbbk}\left[
x_{1},...,x_{n}\right]  /\overline{I}\right)  \leq\sum_{\mathbf{a\in A}%
}\left(  \sum_{i=1}^{n}\nu_{i}\left(  a_{i}\right)  -n+1\right)  .\tag{2.8}%
\end{equation}
By Lemma 2.2 follows that there is a polynomial $g=x_{1}+c_{2}x_{2}+...+c_{n}x_{2}\in
\Bbbk\left[  x_{1},...,x_{n}\right]  ,$ with non zero $c_{2},...,c_{n}\in
\Bbbk$, such that its restriction on $A_{1}\times...\times A_{n}$ is an
injective mapping. Then from $\mu_{\overline{g}}\left(  \overline{g}\right)
=0$ follows that $\mu_{\overline{g}}\left(  g\right)  \in\overline{I}$.
Therefore by Lemma 2.5 we obtain that $\partial_{\mathbf{a}}^{\mathbf{m}}\mu_{\overline{g}%
}\left(  g\right)  $\ $=0,~$for all $\mathbf{a\in A},~\mathbf{m\in}\left[
\mathbf{0,}\nu\left(  \mathbf{a}\right)  -\mathbf{1}\right]  $. For the
special case of $\mathbf{m=0}$, it  follows that $\mu_{\overline{g}}\left(
g\left(  \mathbf{a}\right)  \right)  =0$ for any $\mathbf{a\in A}$.
Therefore, there is a $h\in\Bbbk\left[  t\right]  $ such that
\begin{equation}
\mu_{\overline{g}}\left(  t\right)  =h\left(  t\right)  \prod_{c\in g\left(
\mathbf{A}\right)  }\left(  t-c\right)  .\tag{2.9}%
\end{equation}
Since the restriction of $g$ to $\mathbf{A}$ is an injective function, we have that
$\prod_{c\in g\left(  \mathbf{A}\right)  }\left(  t-c\right)  =\prod
_{\mathbf{a\in A}}\left(  t-g\left(  \mathbf{a}\right)  \right)  $. Therefore
$\left(  2.9\right)  $ becomes
\begin{equation}
\mu_{\overline{g}}\left(  t\right)  =h\left(  t\right)  \prod_{\mathbf{a\in
A}}\left(  t-g\left(  \mathbf{a}\right)  \right)  .\tag{2.10}%
\end{equation}
From $\left(  2.7\right)  $ and $\left(  2.10\right)  $, it is concluded that
for every $\mathbf{a\in A}$ there is a natural number $m_{\mathbf{a}}\in$
$\left[  1,\sum_{i=1}^{n}\nu_{i}\left(  a_{i}\right)  -n+1\right]  $ such
that
\begin{equation}
\mu_{\overline{g}}\left(  t\right)  =\prod_{\mathbf{a\in A}}\left(  t-g\left(
\mathbf{a}\right)  \right)  ^{m_{\mathbf{a}}}.\tag{2.11}%
\end{equation}
Suppose that $m_{\mathbf{a}}<\sum_{i=1}^{n}\nu_{i}\left(  a_{i}\right)  -n+1$ for
some $\mathbf{a\in A,}$ then there exists $\mathbf{m\in}\left[  \mathbf{0,}%
\nu\left(  \mathbf{a}\right)  -\mathbf{1}\right]  $, such that $m_{\mathbf{a}%
}=\left\vert \mathbf{m}\right\vert :=m_{1}+...=m_{n}.$ In this case, using
$\left(  2.11\right)  $ and applying the Leibniz rule for functions of several
variables, we get%
\begin{align*}
& \partial_{\mathbf{a}}^{\mathbf{m}}\mu_{\overline{g}}\left(  g\left(
\mathbf{x}\right)  \right)  \ \\
& =\partial_{\mathbf{a}}^{\mathbf{m}}\left(  \left(  g\left(  \mathbf{x}%
\right)  -g\left(  \mathbf{a}\right)  \right)  ^{m_{\mathbf{a}}}%
\prod_{\substack{\mathbf{b\in A}\\\mathbf{b\neq a}}}\left(  g\left(
\mathbf{x}\right)  -g\left(  \mathbf{b}\right)  \right)  ^{m_{\mathbf{b}}%
}\right)  \\
& =\sum_{\mathbf{r}\in\left[  \mathbf{0,m}\right]  }\dbinom{\mathbf{m}%
}{\mathbf{r}}\partial_{\mathbf{a}}^{\mathbf{r}}\left(  g\left(  \mathbf{x}%
\right)  -g\left(  \mathbf{a}\right)  \right)  ^{m_{\mathbf{a}}}%
\partial_{\mathbf{a}}^{\mathbf{m-r}}\prod_{\substack{\mathbf{b\in
A}\\\mathbf{b\neq a}}}\left(  g\left(  \mathbf{x}\right)  -g\left(
\mathbf{b}\right)  \right)  ^{m_{\mathbf{b}}}\\
=  & m_{\mathbf{a}}!c_{2}^{m_{2}}...c_{n}^{m_{n}}\prod_{\substack{\mathbf{b\in
A}\\\mathbf{b\neq a}}}\left(  g\left(  \mathbf{a}\right)  -g\left(
\mathbf{b}\right)  \right)  ^{m_{\mathbf{b}}}.
\end{align*}
Since $g=x_{1}+\sum_{i=2}^{n}c_{i}x_{i}\,,\,c_{i}\in\Bbbk-\left\{  0\right\}
$ is injective on $\mathbf{A,}$  follows that\
\[
m_{\mathbf{a}}!c_{2}^{m_{2}}...c_{n}^{m_{n}}\prod_{\substack{\mathbf{b\in
A}\\\mathbf{b\neq a}}}\left(  g\left(  \mathbf{a}\right)  -g\left(
\mathbf{b}\right)  \right)  ^{m_{\mathbf{b}}}\neq0.
\]
Therefore $\partial_{\mathbf{a}}^{\mathbf{m}}\mu_{\overline{g}}\left(
g\left(  \mathbf{x}\right)  \right)  \neq0,$ implying that $\mu_{\overline{g}%
}\left(  \overline{g}\right)  \neq\overline{0},$ which is not true. Therefore
for all $\mathbf{a\in A}$ the following is valid
\begin{equation}
m_{\mathbf{a}}=\sum_{i=1}^{n}\nu_{i}\left(  a_{i}\right)  -n+1.\tag{2.12}%
\end{equation}
From $\left(  2.11\right)  $ and $\left(  2.12\right)  $ we get
\[
\mu_{\overline{g}}\left(  t\right)  =\prod_{\mathbf{a\in A}}\left(  t-g\left(
\mathbf{a}\right)  \right)  ^{\sum_{i=1}^{n}\nu_{i}\left(  a_{i}\right)
-n+1},
\]
so one has,%
\[
\deg\mu_{\overline{g}}=\sum_{\mathbf{a\in A}}\sum_{i=1}^{n}\left(  \nu
_{i}\left(  a_{i}\right)  -n+1\right)  .
\]
From this and from $\left(  2.8\right)  $ we get $\left(  2.6\right)  .$
Now from $\left(  2.6\right)  ,$ by Remark $2.1$ it follows that
$\overline{\Bbbk}\left[  x_{1},...,x_{n}\right]  /\overline{I}$ has a
primitive element if and only if%
\[
\dim_{\overline{\Bbbk}}\overline{\Bbbk}\left[  x_{1},...,x_{n}\right]
/\overline{I}=\sum_{\mathbf{a\in A}}\sum_{i=1}^{n}\left(  \nu_{i}\left(
a_{i}\right)  -n+1\right)  ,
\]
which, by Lemma 2.4, can be written equivalently as%
\[
\prod_{i=1}^{n}\sum_{a_{i}\in A_{i}}\nu_{i}\left(  a\right)  =\sum
_{\mathbf{a\in A}}\sum_{i=1}^{n}\left(  \nu_{i}\left(  a_{i}\right)
-n+1\right)
\]
or
\[
\prod_{i=1}^{n}\sum_{a_{i}\in A_{i}}\nu_{i}\left(  a\right)  =\sum
_{\mathbf{a\in A}}\sum_{i=1}^{n}\nu_{i}\left(  a_{i}\right)  -\left(
n-1\right)  \left\vert \mathbf{A}\right\vert .
\]
Changing order of summation on the right-hand side, the above expression becomes
\[
\prod_{i=1}^{n}\sum_{a_{i}\in A_{i}}\nu_{i}\left(  a\right)  =\sum_{i=1}%
^{n}\left(  \frac{\left\vert \mathbf{A}\right\vert }{\left\vert A_{i}%
\right\vert }\sum_{a\in A_{i}}\nu_{i}\left(  a\right)  \right)  -\left(
n+1\right)  \left\vert \mathbf{A}\right\vert .
\]
Dividing the above formula by $\left\vert \mathbf{A}\right\vert ,$ yields%
\begin{equation}
\prod_{i=1}^{n}\frac{\sum_{a_{i}\in A_{i}}\nu_{i}\left(  a\right)
}{\left\vert A_{i}\right\vert }=\sum_{i=1}^{n}\left(  \frac{\sum_{a\in A_{i}%
}\nu_{i}\left(  a\right)  }{\left\vert A_{i}\right\vert }\right)  -\left(
n-1\right)  .\tag{2.13}%
\end{equation}
So we came to the conclusion: $\overline{\Bbbk}\left[  x_{1},...,x_{n}\right]
/\overline{I}$ has a primitive element if and only if $\left(  2.13\right)  $
is valid$.$ Note that $\nu_{i}\left(  a\right)  \geq1,$ $i=1,..,n$ for every
$a\in A_{i}$. Therefore, for any $i=1,..,n$ we have $k_{i}:=\frac{\sum
_{a_{i}\in A_{i}}\nu_{i}\left(  a\right)  }{\left\vert A_{i}\right\vert }%
\geq1$. Thus by Lemma 2.3 we conclude that $\left(  2.13\right)  $ is valid if
and only if at most one of the numbers $\frac{\sum_{a_{i}\in A_{i}}\nu
_{i}\left(  a\right)  }{\left\vert A_{i}\right\vert },$~$i=1,...,n$ is greater
than $1.$ In other words $\overline{\Bbbk}\left[  x_{1},...,x_{n}\right]
/\overline{I}$ has a primitive element if and only if at most one of the
polynomials $f_{1},...,f_{n}$ is inseparable, and in this case $\overline
{g}$ is a primitive element of $\overline{\Bbbk}\left[  x_{1},...,x_{n}%
\right]  /\overline{I},$ where $g=x_{1}+c_{2}x_{2}+...+c_{n}x_{2}\in
\Bbbk\left[  x_{1},...,x_{n}\right]  ,$ with non zero $c_{2},...,c_{n}$.
Finaly by Lemma 2.6 we immediately prove the claim of Theorem 1.1.
\end{proof}

\begin{proof}
[Proof of Corollary 1.2]Let $I$ be the ideal of $\Bbbk\left[  x_{1}%
,...,x_{n}\right]  $ generated by $f_{1}:=\mu_{\delta_{1}}\left(
x_{1}\right)  ,...,f_{1}:=\mu_{\delta_{n}}\left(  x_{n}\right)  $. Then by
Theorem 1.1, $\Bbbk\left[  x_{1},...,x_{n}\right]  /I$ has a primitive element
$\overline{g}$, for some $g=x_{1}+c_{2}x_{2}+...+c_{n}x_{n}\in\Bbbk\left[
x_{1},...,x_{n}\right]  ,$ $\prod_{i=2}^{n}c_{i}\neq0$. Therefore, for every
$f\in\Bbbk\left[  x_{1},...,x_{n}\right]  $ there exist $p\in\Bbbk\left[
t\right]  $, and $h_{1},...,h_{n}\in\Bbbk\left[  x_{1},...,x_{n}\right]  $
such that $f=p\left(  g\right)  +h_{1}f_{1}+...+f_{n}h_{n}$, which by setting
$x_{1}=\delta_{1},...,x_{n}=\delta_{n}$ becomes $f\left(  \delta
_{1},...,\delta_{n}\right)  =p\left(  \delta_{1}+c_{2}\delta_{2}%
+...+c_{n}\delta_{n}\right)  .$
\end{proof}

\begin{proof}
[Proof of Corollary 1.3]From $\left(  1.1\right)  $ we have that
\[
\dim_{\Bbbk
}\Bbbk\left[  \delta_{1},\delta_{2}\right]  =\dim_{\Bbbk}\Bbbk\left[x_{1},x_{2}\right]  /\left(  \mu_{\delta_{1}}\left(  x_{1}\right)
,\mu_{\delta_{2}}\left(  x_{2}\right)  \right).
\]
Therefore $\Bbbk\left[
\delta_{1},\delta_{2}\right]  $ and $\Bbbk\left[  x_{1},x_{2}\right]  /\left(
\mu_{\delta_{1}}\left(  x_{1}\right)  ,\mu_{\delta_{2}}\left(  x_{2}\right)
\right)  $ are isomorphic.$~$Since $\mu_{\delta_{1}},\mu_{\delta_{2}}$ are
both inseparable, by Theorem 1.1 follows that $\Bbbk\left[  x_{1}%
,x_{2}\right]  /\left(  \mu_{\delta_{1}}\left(  x_{1}\right)  ,\mu_{\delta
_{2}}\left(  x_{2}\right)  \right)  $ has no primitive element, meaning that
$\Bbbk\left[  \delta_{1},\delta_{2}\right]  $ has no primitive element.
\end{proof}

\section{The Frobenius problem}

\begin{proof}
[Proof of Theorem 1.4]Since a matrix is diagonalizable if and only if its
minimal polynomial is separable, the proof follows immediately by Corollary 1.2.
\end{proof}

\begin{description}
\item[Lemma 3.1] \bigskip\textit{Let }$A$\textit{\ and}$\ B$\textit{\ be
}$n\times n$\textit{\ commuting matrices over a field }$\Bbbk$\textit{\ with
}$\dim_{\Bbbk}\Bbbk\left[  A,B\right]  =n.$\textit{\ If \ the extension\ }%
$\Bbbk\subset\Bbbk\left[  A,B\right]  $\textit{\ has no primitive element,
then there is no matrix }$C$\textit{\ such that }$A$\textit{\ and}%
$\ B$\textit{\ can be expressed as polynomials in }$C$\textit{.}
\end{description}

\begin{proof}
Suppose there is a $n\times n$ matrix $C\notin\Bbbk\left[  A,B\right]  $ such
that $A$ and$\ B$ are polynomials in $C.$ Thus, $\Bbbk\left[  A,B\right]
\subseteq\Bbbk\left[  C\right]  $, and because by Cayley--Hamilton Theorem is
valid $\dim_{\Bbbk}\Bbbk\left[  C\right]  \leq n,$ we conclude that
$\Bbbk\left[  A,B\right]  =\Bbbk\left[  C\right]  .$ Therefore $C\in
\Bbbk\left[  A,B\right]  $ which contradics to $C\notin\Bbbk\left[  A,B\right]  $.
\end{proof}

\begin{example}
Consider \ the matrices $A$ and$\ B$ as in $\left(  1.1\right)  $. Because
$A^{2}=B^{2}=AB=BA=0$, follows that $\Bbbk\left[  A,B\right]  =\left\{
rI+sA+tB:r,s,t\in\Bbbk\right\}  $, and $\dim_{\Bbbk}\Bbbk\left[  A,B\right]
=3.$ Further through [1] there is no matrix $C\in\Bbbk\left[  A,B\right]  $
such that $A$ and$\ B$ can be expressed as polynomials in $C$. Thus by Lemma
3.1 there is no $3\times3$ matrix $C$ such that $A$ and$\ B$ can be expressed
as polynomials in $C$.
\end{example}

\begin{proof}
[Proof of Theorem 1.5]From $\left(  1.3\right)  $ by Corollary $\mathbf{1.3}$ we obtain the first assertion and then by
Lemma 3.1 we immediately obtain the second assertion.
\end{proof}

\begin{description}
\item[Corollary 3.2] Given the composite number $n$ and the natural numbers
$m,k>1$ such that $mk=n.$ \textit{If }$f=x^{m}+a_{m-1}x^{m-1}+...+a_{0}%
\in\Bbbk\left[  x\right]  $\textit{, }$g=y^{k}+b_{n-1}y^{k-1}+...+b_{0}%
\in\Bbbk\left[  y\right]  \,$\textit{\ are both inseparable, then the
following two }$n\times n$\textit{\ block matrices }$\rho\left(  f\right)
$\textit{\ and }$\rho\left(  g\right)  $\textit{\ can not be expressed as
polynomials of any }$n\times n$\textit{\ matrix }$C,$\textit{\ }%
\[
\rho\left(  f\right)  =\left[
\begin{array}
[c]{ccccccc}%
C\left(  f\right)  & \boldsymbol{0} & \boldsymbol{0} & . & . & . &
\boldsymbol{0}\\
\boldsymbol{0} & C\left(  f\right)  & \boldsymbol{0} & . & . & . &
\boldsymbol{0}\\
. & . & . & . &  &  & .\\
. & . & . &  & . &  & .\\
. & . & . &  &  & . & .\\
\boldsymbol{0} & \boldsymbol{0} & \boldsymbol{0} & . & . & . & \boldsymbol{0}%
\\
\boldsymbol{0} & \boldsymbol{0} & \boldsymbol{0} & . & . & . & C\left(
f\right)
\end{array}
\right]
\]
$~$
\end{description}

\textit{and}
\[
\rho\left(  g\right)  =\left[
\begin{array}
[c]{ccccccc}%
\boldsymbol{0} & \boldsymbol{I} & \boldsymbol{0} & . & . & . & \boldsymbol{0}%
\\
\boldsymbol{0} & \boldsymbol{0} & \boldsymbol{I} & . & . & . & \boldsymbol{0}%
\\
. & . & . & . &  &  & .\\
. & . & . &  & . &  & .\\
. & . & . &  &  & . & .\\
\boldsymbol{0} & \boldsymbol{0} & \boldsymbol{0} & . & . & . & \boldsymbol{I}%
\\
-b_{0}\boldsymbol{I} & -b_{1}\boldsymbol{I} & -b_{2}\boldsymbol{I} & . & . &
. & -b_{k}\boldsymbol{I}%
\end{array}
\right]  ,
\]

\textit{where }$\boldsymbol{0}$\textit{\ and }$\boldsymbol{I}$\textit{\ are
the zero and unit }$m\times m$\textit{\ matrices, respectively, and }$C\left(
f\right)  $\textit{\ is the companion matrix of the polynomial }%
$f$\textit{\ given by }%
\[
C\left(  f\right)  =\left[
\begin{array}
[c]{ccccccc}%
0 & 1 & 0 & . & . & . & 0\\
0 & 0 & 1 & . & . & . & 0\\
. & . & . & . &  &  & .\\
. & . & . &  & . &  & .\\
. & . & . &  &  & . & .\\
0 & 0 & 0 & . & . & . & 1\\
-a_{0} & -a_{1} & -a_{2} & . & . & . & -a_{m}%
\end{array}
\right]  .
\]

\textit{Further the following holds:} \textit{Let }$\mathbf{M}$\textit{ be the
set of pairs }$\left(  f,g\right)  $\textit{ of all inseparable polynomials
}$f$\textit{ and }$g~$\textit{satisfying the condition }$\deg f\deg
g=n$\textit{.~Then the set of all matrix pairs }$\left(  A,B\right)  $\textit{
satisfying the conditions }$\left(  1.3\right)  $\textit{ of Theorem 1.5 is
given by }$\left\{  \left(  S\rho\left(  f\right)  S^{-1},S\rho\left(
g\right)  S^{-1}\right)  :S\in\operatorname*{GL}_{n}\left(  \Bbbk\right)
,~\left(  f,g\right)  \in\mathbf{M}\right\}  .$

\begin{proof}
A basis of $\Bbbk\left[  x,y\right]  /\left(  f,g\right)  $ as vector space
is
\[
B=\left\{  \widetilde{x}^{i}\widetilde{y}^{j}%
:i=0,1,...,m-1;j=0,1,...,k-1~\right\}  .
\]
By ordering the basis $B$ as follows,
\begin{equation}
B=\left\{  1,\widetilde{x},...,\widetilde{x}^{m-1},\widetilde{y},\widetilde
{x}y,...,\widetilde{x}^{m-1}\widetilde{y},...,\widetilde{y}^{k-1}%
,\widetilde{x}\widetilde{y}^{k-1},...,\widetilde{x}^{m-1}\widetilde{y}%
^{k-1}\right\}  ,\tag{3.1}%
\end{equation}
it is easy to verify, that the regular matrix representations of
$\widetilde{x}$ and $\widetilde{y}$ with respect to basis $B$, as ordered in
$\left(  3.1\right)  ,$ equal to the matrices $\rho\left(  f\right)
$\textit{\ and }$\rho\left(  g\right)  $ respectively. Now the first assertion
follows by Theorem 1.5.
Finally, the second claim follows directly from the first and we will omit the details.
\end{proof}

\begin{example}
Let $f=\left(  x-1\right)  ^{2}\left(  x+2\right)  =\allowbreak x^{3}-3x+2$
and $g=\left(  x+1\right)  ^{2}=\allowbreak x^{2}+2x+1$, then by Corollary 3.2
the matrices $\rho\left(  f\right)  =\left[
\begin{array}
[c]{cccccc}%
0 & 1 & 0 & 0 & 0 & 0\\
0 & 0 & 1 & 0 & 0 & 0\\
-2 & 3 & 0 & 0 & 0 & 0\\
0 & 0 & 0 & 0 & 1 & 0\\
0 & 0 & 0 & 0 & 0 & 1\\
0 & 0 & 0 & -2 & 3 & 0
\end{array}
\right]  $ and $\rho\left(  g\right)  =\left[
\begin{array}
[c]{cccccc}%
0 & 0 & 0 & 1 & 0 & 0\\
0 & 0 & 0 & 0 & 1 & 0\\
0 & 0 & 0 & 0 & 0 & 1\\
-1 & 0 & 0 & -2 & 0 & 0\\
0 & -1 & 0 & 0 & -2 & 0\\
0 & 0 & -1 & 0 & 0 & -2
\end{array}
\right]  $ $\ $can not be expressed as polynomials in any one matrix $C.$
\end{example}


\begin{thebibliography}{9}                                                                                                %


\bibitem {[1]}\bigskip Frobenius, F. G., Ueber vertauschbare Matrizen,
\textit{Sitzungsberichte der Akad. der Wiss. zu Berlin}, 1896. Persistent
Link: http://dx.doi.org/10.3931/e-rara-18867

\bibitem {[2]}\bigskip M. Gerstenhaber, On Dominance and Varieties of
Commuting Matrices, \textit{Annals of Mathematics} , Vol. 73, No. 2, pp.
324-548, March 1961.

\bibitem {[3]}Horn, R.A, \& Johnson, C.R, Matrix Analysis,\textit{ Cambridge
University Press},1985

\bibitem {[4]}Kechriniotis, A. I., Delibasis, K. K., Oikonomou, I. P., \&
Tsigaridas, G. N. (2023). Classical multivariate Hermite coordinate
interpolation on n-dimensional grids, \textit{ArXiv preprint arXiv:2301.01833}.

\bibitem {[5]}K.Shoda, Ueber mit einer matrix vertauschbare matrizen,
\textit{Mathematische Zeitschrift}, vol. 29, pp. 696-712, 1929.
\end{thebibliography}
\end{document}